\documentclass[11pt,a4j]{article}

\usepackage{lineno,hyperref}
 \usepackage{amsmath,amssymb,amsfonts}
 \usepackage{amsthm}

\newtheorem{theorem}{Theorem}[section]

\newtheorem{prop}[theorem]{Proposition}
\newtheorem{rem}{Remark}
\newtheorem{example}{Example}
\newtheorem{defn}{Definition}

\title{The deformations of symplectic structures by moment maps}
\author{Tomoya Nakamura  \thanks{Department of Mathematics, Waseda University, \textup{3-4-1}, Okubo, Shinjuku-ku, Tokyo, Japan}
 \\ email: \href{mailto:x-haze@ruri.waseda.jp}{x-haze@ruri.waseda.jp}}
\date{\today}

\begin{document}

  \maketitle
  
  \begin{abstract}
  We study deformations of symplectic structures on a smooth manifold $M$ via the quasi-Poisson theory.
By a fact, we can deform a given symplectic structure $\omega $ to a new symplectic structure $\omega _t$ parametrized by some element $t$ in $\Lambda^2\mathfrak{g}$, where $\mathfrak{g}$ is the Lie algebra of a Lie group $G$. Moreover, we can get a lot of concrete examples for the deformations of symplectic structures on the complex projective space and the complex Grassmannian.
  \end{abstract}

 \section{Introduction}\label{sec:1}
In the context of symplectic geometry, deformation-equivalence assumptions and conditions are often appeared, for example, in the statement of Moser's theorem \cite{Si} and Donaldson's four-six conjecture \cite{Sa}. However, it seems that a method of constructing deformation-equivalent symplectic structures specifically is not well known. In this paper, we construct a method of producing new symplectic structures deformation-equivalent to a given symplectic structure. %Two symplectic structures $\omega _0$ and $\omega _1$ are {\it deformation-equivalent} if there is a smooth family $\omega _t$ of symplectic structures joining $\omega _0$ and $\omega _1$. 
Our approach to deformations of symplectic structures is to use quasi-Poisson theory which was introduced by Alekseev and Kosmann-Schwarzbach \cite{AK}, and this approach is carried out by using the fact that a moment map for a symplectic-Hamiltonian action $\sigma $ is also a moment map for a quasi-Poisson action $\sigma $. The former moment map satisfies conditions for only one symplectic structure, whereas the latter does conditions for a family of quasi-Poisson structures parametrized by elements in $\Lambda ^2\mathfrak{g}$. From here we call these elements {\it twists}.  Regarding the quasi-Poisson structure induced by a symplectic structure as that with twist $0$, which is denoted by $\pi _0$,  we can find different quasi-Poisson structures $\pi _t$ which induce symplectic structures $\omega _t$ by the choice of "good" twists $t$. The quasi-Poisson structure inducing a symplectic structure must be a nondegenerate Poisson structure. We describe the conditions for the quasi-Poisson structure with a twist $t$ to be a nondegenerate Poisson structure. Our method of using the family of quasi-Poisson structures is one of interesting geometry frameworks (See \cite{AK}).

From here, we explain briefly the difference among moment maps for symplectic- and quasi-Poisson-Hamiltonian actions, and equivariant moment maps for Poisson actions on a smooth manifold (In Poisson geometry, {\it non}-equivariant moment maps for Poisson action can be defined \cite{L}, \cite{L2}.).\\
(I) Symplectic-Hamiltonian actions

In symplectic geometry, a moment map $\mu :M\rightarrow \mathfrak{g}^*$ for a symplectic action $\sigma $ of a Lie group $G$ on a symplectic manifold $(M,\omega )$ is defined with two conditions: one is for the symplectic structure $\omega $,
$$d\mu ^X=\iota _{X_M}\omega \ (X\in \mathfrak{g}).$$
Here $\mu ^X(p):=<\mu (p),X>$ and $X_M$ is a vector field on $M$ defined by
\begin{equation}\label{vector field}
X_{M,p}:=\left. \frac{d}{dt}\sigma _{\exp tX}(p)\right|_{t=0}
\end{equation}
for $p$ in $ M$. The other is the $G$-equivariance condition with respect to the action $\sigma $ on $M$ and the coadjoint action $\mathrm{Ad}^*$ on $\mathfrak{g}^*$,
$$\mu \circ \sigma _g=\mathrm{Ad}^*_g\circ \mu $$
for all $g$ in $G$. In this paper, we call symplectic actions with moment maps {\it symplectic-Hamiltonian actions
} to distinguish it from other actions with moment maps.\\
(II) Poisson actions with equivariant moment maps

A Poisson Lie group, which was introduced by Drinfel'd \cite{D}, is a Lie group with a Poisson structure $\pi $ compatible with the group structure. Namely, the structure $\pi $ satisfies
\begin{eqnarray}\label{multiplicative}
\pi _{gh}=L_{g*}\pi _h+R_{h*}\pi _g
\end{eqnarray}
for any $g$ and $h$ in $G$, where $L_g$ and $R_h$ are the left and right translations in $G$ by $g$ and $h$, respectively. Such a structure is called {\it multiplicative}. Then the simply connected Lie group $G^*$ called the dual Poisson Lie group is obtained uniquely from a Poisson Lie group $(G,\pi )$ and a local action of $G$ on $G^*$ is defined naturally. We call a multiplicative Poisson structure $\pi $ on $G$ {\it complete} if the action is global. Then $(G,\pi )$ is called a {\it complete Poisson Lie group}. An equivariant moment map $\mu :M\rightarrow G^*$ for a Poisson action $\sigma $ of a complete Poisson Lie group $(G,\pi )$ on a Poisson manifold $(M,\pi )$ is a generalization of a moment map for a symplectic-Hamiltonian action on a symplectic manifold, which was given by Lu in \cite{L}.\\
(III) Quasi-Poisson-Hamiltonian actions

Quasi-Poisson theory, which was originated with \cite{AK} by Alekseev and Kosmann-Schwarzbach, is a generalization of Poisson theory with Poison actions. In quasi-Poisson geometry, quasi-triples $(D,G,\mathfrak{h})$ and its infinitesimal version, Manin quasi-triples $(\mathfrak{d},\mathfrak{g},\mathfrak{h})$, play important roles. A quasi-triple $(D,G,\mathfrak{h})$ defines a quasi-Poisson Lie group $G_D^\mathfrak{h}$ and we can obtain the notion of a quasi-Poisson action of such a quasi-Poisson Lie group $G_D^\mathfrak{h}$. A moment map $\mu $ for the action is a map from $M$ into $D/G$ and satisfies a condition not for one quasi-Poisson structure but for a family of quasi-Poisson structures parametrized by elements in $\Lambda ^2\mathfrak{g}$. An equivariant moment map for a Poisson action in (II) is an example of a moment map for a quasi-Poisson-Hamiltonian action. In this paper, we use the moment map theory for quasi-Poisson actions to deform symplectic structures on a smooth manifold.

This paper is constructed as follows. It is contents of Section 2 to review the moment map theory for quasi-Poisson actions. In Section 3, we describe a deformation method of symplectic structures on a smooth manifold via the quasi-Poisson theory. This method is the subject in this paper. %A twist $t$ in $\Lambda ^2\mathfrak{g}$ deforms an isotropic complement $\mathfrak{g}^*$ of $\mathfrak{g}$ in $\mathfrak{g}\oplus \mathfrak{g}^*$ to $\mathfrak{g}^*_t$. By observation of the moment map condition for each space $\mathfrak{g}^*_t$, we can decide which twist gives a deformation between symplectic structures. 
Theorem \ref{abstract main} %gives conditions for a twist to deform the nondegenerate (quasi-)Poisson structure induced by a symplectic structure to a Poisson structure and to a nondegenerate quasi-Poisson structure. Since any nondegenerate Poisson structure induces a symplectic structure, the theorem also 
gives a condition for a twist to deform a symplectic structure to a new one. In addition, Theorem \ref{main} gives a sufficient condition for a twist to satisfy the assumption of Theorem \ref{abstract main}.  In Section 4, we introduce concrete examples for deformations of symplectic structures. %First, we deform the Fubini-Study form on the complex projective space $\mathbb{C}\mathrm{P}^k$ with an action of $\mathrm{SU}(k+1)$, where $k=1,2$. Second example is a deformation of the Fubini-Study form on $\mathbb{C}\mathrm{P}^n$with an action of the torus group $\mathbb{T}^n$, which remains the symplectic-Hamiltonian action with the same moment map after the deformation. The last one is the complex Grassmannian $\mathrm{Gr}(r;\mathbb{C}^n)$ with the Kirillov-Kostant form. We show that it is deformed by some twist.
% the twist $\lambda t$, where
%$$t=\frac{1}{4n}\sum _{1\leq i<j\leq n}^{}X_{ij}\wedge Y_{ij}\in \Lambda ^2\mathfrak{su}(n)$$
%
%and $\lambda $ is a sufficiently small real number. Here $X_{ij}$ and $Y_{ij}$ are elements in some basis of $\mathfrak{su}(n)$.
We give deformations of the Fubini-Study and the Kirillov-Kostant forms on $\mathbb{C}\mathrm{P}^n$ and the complex Grassmannian, respectively. 

\section{Moment maps for quasi-Poisson actions on quasi-Poisson manifolds}\label{quasi-Poisson}

In this section, we shall recall the quasi-Poisson theory \cite{AK}. We start with the definition of quasi-Poisson Lie groups, which is a generalization of Poisson Lie groups.

\begin{defn}\label{q-PL}
Let $G$ be a Lie group with the Lie algebra $\mathfrak{g}$. Then a pair $(\pi ,\varphi )$ is a {\it quasi-Poisson structure} on $G$ if a multiplicative $2$-vector field $\pi $ on $G$ and an element
$\varphi $ of $\Lambda ^3\mathfrak{g}$ satisfy
\begin{eqnarray}
\frac{1}{2}\left[\pi ,\pi \right]&=&\varphi ^R-\varphi ^L,
 \label{q-PL-1}\\ 
\left[\pi ,\varphi^L\right]&=&\left[\pi ,\varphi ^R\right]=0, \label{q-PL-2}
\end{eqnarray}
where the bracket $[\cdot ,\cdot ]$ is the Schouten bracket on the multi-vector fields on $G$, and $\varphi ^L$ and $\varphi ^R$ denote the left and right invariant $2$-vector fields on $G$ with value $\varphi $ at $e$ respectively. A triple $(G, \pi ,\varphi )$ is called a {\it quasi-Poisson Lie group}.
\end{defn}

\begin{rem}
In a quasi-Poisson structure $(\pi ,\varphi )$ on $G$, the 2-vector field $\pi $ is a multiplicative Poisson structure if $\varphi =0$. Namely, $(G,\pi )$ is a Poisson Lie group. 
\end{rem}

We use a "quasi-triple" to obtain a quasi-Poisson Lie group. To define a quasi-triple, we describe its infinitesimal version, a Manin quasi-triple. 

\begin{defn}
Let $\mathfrak{d}$ be a $2n$-dimensional Lie algebra with an invariant nondegenerate symmetric bilinear form of signature $(n,n)$, which is denoted by $(\cdot |\cdot )$. Let $\mathfrak{g}$ be an $n$-dimensional Lie subalgebra of $\mathfrak{d}$ and $\mathfrak{h}$ be an $n$-dimensional vector subspace of $\mathfrak{d}$. 
Then a triple $(\mathfrak{d},\mathfrak{g},\mathfrak{h})$ is a {\it Manin quasi-triple} if $\mathfrak{g}$ is a maximal isotropic subspace with respect to $(\cdot |\cdot )$ and $\mathfrak{h}$ is an isotropic complement subspace of $\mathfrak{g}$ in $\mathfrak{d}$.
\end{defn}

\begin{rem}
For a given Lie algebra $\mathfrak{d}$ and a Lie subalgebra $\mathfrak{g}$ of $\mathfrak{d}$, a choice of an isotropic complement subspace $\mathfrak{h}$ of $\mathfrak{g}$ in $\mathfrak{d}$ is not unique.
\end{rem}

A Manin quasi-triple $(\mathfrak{d},\mathfrak{g},\mathfrak{h})$ defines the decomposition $\mathfrak{d}=\mathfrak{g}\oplus \mathfrak{h}$. Then the linear isomorphism
\begin{eqnarray}\label{linear isomorphism}
\ j:\mathfrak{g}^* \rightarrow \mathfrak{h},\ (j(\xi )|x):=<\xi ,x>\ (\xi \in \mathfrak{g}^*,x\in \mathfrak{g})
\end{eqnarray}
is determined by the decomposition. We denote the projections from $\mathfrak{d}=\mathfrak{g}\oplus \mathfrak{h}$ to $\mathfrak{g}$ and $\mathfrak{h}$ by $\ p_\mathfrak{g}$ and $p_\mathfrak{h}$ respectively. We introduce an element $\varphi _\mathfrak{h}$ in $ \Lambda ^3\mathfrak{g}$ which is defined by the map from $\Lambda ^2\mathfrak{g}^*$ to $\mathfrak{g}$, denoted by the same letter,
\begin{eqnarray}\label{varphi}
\varphi _\mathfrak{h}(\xi ,\eta )=p_\mathfrak{g}(\left[j(\xi ),j(\eta )\right]),
\end{eqnarray}
for any $\xi ,\eta $ in $ \mathfrak{g}^*$. We define the linear map $F_\mathfrak{h}:\mathfrak{g}\rightarrow \Lambda ^2\mathfrak{g}$ by setting
\begin{eqnarray}\label{cobracket dual}
F_\mathfrak{h}^*(\xi ,\eta )=j^{-1}(p_\mathfrak{h}(\left[j(\xi ),j(\eta )\right]))
\end{eqnarray}
for any $\xi ,\eta $ in $ \mathfrak{g}^*$, where $F_\mathfrak{h}^*:\Lambda ^2\mathfrak{g}^*\rightarrow \mathfrak{g}^*$ is the dual map of $F_\mathfrak{h}$. These elements will be used later to define a quasi-Poisson structure and a quasi-Poisson action respectively.

Next we define a quasi-triple $(D,G,\mathfrak{h})$ and construct a quasi-Poisson structure on $G$ using $(D,G,\mathfrak{h})$. 

\begin{defn}
Let $D$ be a connected Lie group with a bi-invariant scalar product with the Lie algebra $\mathfrak{d}$ and $G$ be a connected closed Lie subgroup of $D$ with the Lie algebra $\mathfrak{g}$. Let $\mathfrak{h}$ be a vector subspace of $\mathfrak{d}$.  Then a triple $(D,G,\mathfrak{h})$ is a {\it quasi-triple} if $(\mathfrak{d},\mathfrak{g},\mathfrak{h})$ is a Manin quasi-triple.
\end{defn}

A method of constructing a quasi-Poisson structure by a quasi-triple is as follows. Let $(D,G,\mathfrak{h})$ be a quasi-triple with a Manin quasi-triple $(\mathfrak{d},\mathfrak{g},\mathfrak{h})$. Using the inverse $j^{-1}:\mathfrak{h}\rightarrow \mathfrak{g}^*$ of the linear isomorphism (\ref{linear isomorphism}), we identify $\mathfrak{d}$ with $\mathfrak{g}\oplus \mathfrak{g}^*$. Consider the map
$$r_\mathfrak{h}:\mathfrak{d}^*\rightarrow \mathfrak{d},\ \xi +X\mapsto \xi,$$
for any $\xi $ in $ \mathfrak{g}^*$ and $X$ in $ \mathfrak{g}$. This map defines an element $r_\mathfrak{h}\in \mathfrak{d}\otimes \mathfrak{d}$ which we denote by the same letter. We set
$$\pi _D^\mathfrak{h}:=r_\mathfrak{h}^L-r_\mathfrak{h}^R,$$
where $r_\mathfrak{h}^L$ and $r_\mathfrak{h}^R$ is denoted as the left and right invariant $2$-tensors on $D$ with value $r_\mathfrak{h}$ at the identity element $e$ in $D$ respectively, and we can see that it is a multiplicative $2$-vector field on $D$. Furthermore, the $2$-vector field $\pi _D^\mathfrak{h}$ and the element $\varphi _\mathfrak{h}$ defined by (\ref{varphi}) satisfy (\ref{q-PL-1}) and (\ref{q-PL-2}). We set
\begin{eqnarray}\label{q-PofG}
\pi _{G,g}^\mathfrak{h}:=\pi _{D,g}^\mathfrak{h}
\end{eqnarray} 
for any $g$ in $ G$. Then we can see that $\pi _G^\mathfrak{h}$ is well-defined and that $\pi _G^\mathfrak{h}$ is a multiplicative $2$-vector field on $G$. Moreover, $\pi _G^\mathfrak{h}$ and $\varphi _\mathfrak{h}$ satisfy (\ref{q-PL-1}) and (\ref{q-PL-2}). Therefore $(G,\pi _G^\mathfrak{h}, \varphi _\mathfrak{h})$ is a quasi-Poisson Lie group. We sometimes denote a Lie group with such a structure by $G_D^\mathfrak{h}$. 

From here, we consider only connected quasi-Poisson Lie group $G_D^\mathfrak{h}$ defined as above by a quasi-triple $(D,G,\mathfrak{h})$. For a smooth manifold $M$ with a $2$-vector field $\pi _M$, a quasi-Poisson action is defined as follows. It is a generalization of Poisson actions of connected Poisson Lie groups \cite{LW}.

\begin{defn}\label{q-P-mfd}
Let $G_D^\mathfrak{h}$ be a connected quasi-Poisson Lie group acting on a smooth manifold $M$ with a $2$-vector field $\pi _M$. The action $\sigma $ of $G$ on $M$ is a {\it quasi-Poisson action} if for each $X$ in $ \mathfrak{g}$,
\begin{eqnarray}
\frac{1}{2}\left[\pi _M,\pi _M\right]&=&(\varphi _\mathfrak{h})_M, \label{q-P-a-1}\\
\mathfrak{L}_{X_M}\pi _M&=&F_\mathfrak{h}(X)_M, \label{q-P-a-2}
\end{eqnarray}
where $x_M$ is a fundamental multi-vector field for any $x$ in $ \wedge ^*\mathfrak{g}$. Here $F_\mathfrak{h}$ is the dual of the map (\ref{cobracket dual}). Then a $2$-vector field $\pi _M$ is called a {\it quasi-Poisson $G_D^\mathfrak{h}$-structure} on $M$ and $(M,\pi _M)$ is called a {\it quasi-Poisson $G_D^\mathfrak{h}$-manifold}.
\end{defn}

\begin{rem}
A quasi-Poisson Lie group $G_D^\mathfrak{h}$ with the natural left action is not a quasi-Poisson $G_D^\mathfrak{h}$-manifold. In fact, $(\varphi _\mathfrak{h})_G=\varphi _\mathfrak{h}^R$.
\end{rem}

Finally we define a moment map for a quasi-Poisson action to carry out the deformation of symplectic structures using the moment map theory for quasi-Poisson actions in Section \ref{Main Result}. We need some preliminaries to define a moment map. For any quasi-triple $(D,G,\mathfrak{h})$, since $G$ is a closed subgroup of $D$, the quotient space $D/G$ is a smooth manifold, which is the range of moment maps. The action of $D$ on itself by left multiplication induces an action of $D$ on $D/G$. We call it {\it dressing action} of $D$ on $D/G$ and denote the corresponding infinitesimal action by $X\mapsto X_{D\!/\!G}$ for $X$ in $ \mathfrak{d}$. Let $p_{D\!/\!G}:D\rightarrow D/G$ be the natural projection. Then
$$\pi _{D\!/\!G}^\mathfrak{h}:=p_{D\!/\!G*}\pi _D^\mathfrak{h}$$
is a $2$-vector field on $D/G$. We consider the dressing action on $D/G$ restricted to $G$, and can see that $\pi _{D\!/\!G}^\mathfrak{h}$ satisfies (\ref{q-P-a-1}) and (\ref{q-P-a-2}). Therefore $(D/G,\pi _{D\!/\!G}^\mathfrak{h})$ is a quasi-Poisson $G_D^\mathfrak{h}$-manifold. The following definition is one of the important notions to define moment maps.

\begin{defn}\label{admissible}
An isotropic complement $\mathfrak{h}$ of $\mathfrak{g}$ in $\mathfrak{d}$ is called {\it admissible} at a point $s$ in $D/G$ if the infinitesimal dressing action restricted to $\mathfrak{h}$ defines an isomorphism from $\mathfrak{h}$ onto $T_s(D/G)$, that is, the map $\mathfrak{h}\rightarrow T_s(D/G),\ \xi \mapsto \xi _{D\!/\!G,s}$ is an isomorphism. A quasi-triple $(D,G,\mathfrak{h})$ is {\it complete} if $\mathfrak{h}$ is admissible everywhere on $D/G$.
\end{defn}

\noindent
It is clear that any isotropic complement $\mathfrak{h}$ of $\mathfrak{g}$ is admissible at $eG$ in $D/G$. If the complement $\mathfrak{h}$ is admissible at a point $s$ in $D/G$, then it is also admissible on some open neighborhood $U$ of $s$. For a quasi-triple $(D,G,\mathfrak{h})$, we assume that $\mathfrak{h}$ is admissible on an open subset $U$ of $D/G$. Then for any $X$ in $ \mathfrak{g}$, we define the $1$-form $\hat{X}_\mathfrak{h}$ on $U$ by
\begin{eqnarray}\label{1-form}
<\hat{X}_\mathfrak{h},\xi _{D\!/\!G}>=(X|\ \xi )
\end{eqnarray}
for any $\xi $ in $ \mathfrak{h}$. If a quasi-triple $(D,G,\mathfrak{h})$ is complete, then $\hat{X}_\mathfrak{h}$ is a global $1$-form on $D/G$. Next we define a twist between isotropic complement subspaces $\mathfrak{h}$ and $\mathfrak{h}'$ of $\mathfrak{g}$ in $\mathfrak{d}$. Twists also play an important role in the moment map theory for quasi-Poisson actions. Let $j$ and $j'$ be the linear isomorphism (\ref{linear isomorphism}) defined by Manin quasi-triples $(\mathfrak{d},\mathfrak{g},\mathfrak{h})$ and $(\mathfrak{d},\mathfrak{g},\mathfrak{h}')$ respectively. Consider the map
$$t:=j'-j: \mathfrak{g}^*\rightarrow \mathfrak{d}.$$
It is easy to show that $t$ takes values in $\mathfrak{g}$ and that it is anti-symmetric, so that the map $t$ defines an element $t$ in $ \Lambda ^2\mathfrak{g}$ which we denote by the same letter. The element $t$ is called the {\it twist} from $\mathfrak{h}$ to $\mathfrak{h}'$. Fix a quasi-triple $(D,G,\mathfrak{h})$. Let $\mathfrak{h}_t$ be an isotropic complement of $\mathfrak{g}$ with a twist $t$ from $\mathfrak{h}$. Then we can represent the elements $\varphi _{\mathfrak{h}_t}, F_{\mathfrak{h}_t}$ and $\pi _G^{\mathfrak{h}_t}$ defined by a quasi-triple $(D,G,\mathfrak{h}_t)$ as follows:
\begin{eqnarray}
\varphi _{\mathfrak{h}_t}&=&\varphi _\mathfrak{h}+\frac{1}{2}[t,t]+\varphi _t,\label{phi-t}\\
F_{\mathfrak{h}_t}&=&F_\mathfrak{h}+F_t,\label{F-t}\\
\pi _G^{\mathfrak{h}_t}&=&\pi _G^\mathfrak{h}+t^L-t^R,\label{pi_G-t}
\end{eqnarray} 
where $[t,t]:=[t^L,t^L]_e$, $\varphi _t(\xi ):=\overline{\mathrm{ad}_\xi ^*t}$ and $F_t(X):=\mathrm{ad}_Xt$. Here $\mathrm{ad}$ denotes the adjoint action of $\mathfrak{g}$ on $\Lambda ^2\mathfrak{g}$ and $\overline{\mathrm{ad}_\xi ^*t}$ denotes the projection of $\mathrm{ad}_\xi ^*t$ onto $\Lambda ^2\mathfrak{g}\subset \Lambda ^2\mathfrak{d}$, where $\mathfrak{d}^*$ including $\mathfrak{g}^*$ acts on $\Lambda ^2\mathfrak{d}$ by the coadjoint action. Let $\{e_i\}$ be a basis on $\mathfrak{g}$ and $\{\varepsilon ^i\}$ be the basis on $\mathfrak{h}$ identified with the dual basis of $\{e_i\}$ on $\mathfrak{g}^*$ by $j^{-1}$. Then the basis $\{\varepsilon_t^i\}$ on $\mathfrak{h}_t$ identified with the dual basis of $\{e_i\}$ on $\mathfrak{g}^*$ by ${j'}^{-1}$ can be written by
\begin{eqnarray}\label{basis}
\varepsilon_t^i=\varepsilon ^i+t^{ij}e_j,
\end{eqnarray}
where $t=\frac{1}{2}t^{ij}e_i\wedge e_j$. Moreover components of $\varphi _t$ with respect to the basis $\{\varepsilon ^i\}$ are written as
\begin{eqnarray}\label{components}
\varphi _t^{ijk}=(F_\mathfrak{h})_l^{jk}t^{il}-(F_\mathfrak{h})_l^{ik}t^{jl}.
\end{eqnarray}
This indication is useful later. Let $(M,\pi _M^\mathfrak{h})$ be a quasi-Poisson $G_D^\mathfrak{h}$-manifold. We set that $\pi _M^{\mathfrak{h}_t}:=\pi _M^\mathfrak{h}-t_M$. Then it follows that $(M,\pi _M^{\mathfrak{h}_t})$ is a quasi-Poisson $G_D^{\mathfrak{h}_t}$-manifold. Now we define moment maps for quasi-Poisson actions.

\begin{defn}
Let $G_D^\mathfrak{h}$ be a connected quasi-Poisson Lie group defined by a quasi-triple $(D,G,\mathfrak{h})$ and $(M,\pi _M^\mathfrak{h})$ be a quasi-Poisson $G_D^\mathfrak{h}$-manifold. Then a map $\mu :M\rightarrow D/G$ which is equivariant with respect to the $G$-action on $M$ and the dressing action on $D/G$ is a {\it moment map} for the quasi-Poisson action of $G_D^\mathfrak{h}$ on $(M,\pi _M^\mathfrak{h})$ if for any open subset $\Omega \subset M$ and any isotropic complement $\mathfrak{h}'$ admissible on $\mu (\Omega )$,
\begin{eqnarray}\label{q-P-mm}
(\pi _M^{\mathfrak{h}'})^\sharp (\mu ^*(\hat{X}_{\mathfrak{h}'}))=X_M
\end{eqnarray}
on $\Omega$ for any $X$ in $ \mathfrak{g}$. Here $<(\pi _M^{\mathfrak{h}'})^\sharp (\alpha ),\beta >:=\pi _M^{\mathfrak{h}'}(\alpha ,\beta )$. We call a quasi-Poisson action with a moment map a {\it quasi-Poisson-Hamiltonian action}.
\end{defn}

Actually we need not impose the equation (\ref{q-P-mm}) on all admissible complements because we have the following proposition.
 
\begin{prop}[\cite{AK}]\label{sufficient}
Let $\mathfrak{h}$ and $\mathfrak{h}'$ be two complements admissible at a point $s$ in $D/G$, and $p$ in $ M$ be such that $\mu (p)=s$. Then, at the point $p$, conditions (\ref{q-P-mm}) for $\mathfrak{h}$ and $\mathfrak{h}'$ are equivalent, namely 
$$(\pi _M^\mathfrak{h})^\sharp (\mu ^*(\hat{X}_\mathfrak{h}))_p=(\pi _M^{\mathfrak{h}'})^\sharp (\mu ^*(\hat{X}_{\mathfrak{h}'}))_p.$$
\end{prop}

For a quasi-Poisson manifold with a quasi-Poisson-Hamiltonian action, the following theorem holds.

\begin{theorem}[\cite{AK}]\label{generalized foliation}
 Let $(M,\pi _M^\mathfrak{h})$ be a quasi-Poisson manifold on which a quasi-Poisson Lie group $G_D^\mathfrak{h}$ defined by a quasi-triple $(D,G,\mathfrak{h})$ acts by a quasi-Poisson-Hamiltonian action $\sigma $ . For any $p$ in $ M$, if both $\mathfrak{h}'$ and $\mathfrak{h}''$ are admissible at $\mu (p)$ in $D/G$, then 
$$\mathrm{Im}(\pi _M^{\mathfrak{h}'})_p^\sharp =\mathrm{Im}(\pi _M^{\mathfrak{h}''})_p^\sharp ,$$
where $\mu $ is a moment map for $\sigma $.
\end{theorem}

Now we show important examples for quasi-Poisson-Hamiltonian actions.

\begin{example}[Poisson manifolds \cite{AK},\cite{BC},\cite{LW}]\label{Ex.Poisson}
Let $(M,\pi )$ be a Poisson manifold on which a connected Poisson Lie group $(G,\pi _G)$ acts by a Poisson action $\sigma $. Then $(M,\pi )$ is a quasi-Poisson $(G,\pi _G,0)$-manifold and $\sigma $ is a quasi-Poisson action on $(M,\pi )$. In fact, the Manin triple $(\mathfrak{g}\oplus \mathfrak{g}^*,\mathfrak{g},\mathfrak{g}^*)$ corresponding to $(G,\pi _G)$ is a Manin quasi-triple and the multiplicative $2$-vector field $\pi _G$ on $G$ coincides with the $2$-vector field $\pi _G^{\mathfrak{g}^*}$ defined by the corresponding quasi-triple $(D,G,\mathfrak{g}^*)$. Since $[\pi ,\pi ]=0$ and the Poisson action $\sigma $ satisfies
\begin{eqnarray}
\mathfrak{L}_{X_M}\pi =F_{\mathfrak{g}^*}(X)_M
\end{eqnarray}
for any $X$ in $\mathfrak{g}$, the action $\sigma $ is a quasi-Poisson action by Definition \ref{q-P-mfd}. Here the dual of $F_{\mathfrak{g}^*}$ coincides with the bracket on $\mathfrak{g}^*$ defined by $(G,\pi _G)$.

We assume that $\pi _G$ is complete and that there exists a $G$-equivariant moment map $\mu :M\rightarrow G^*$ for the Poisson action $\sigma $, where $G^*$ is the dual Poisson Lie group of $(G,\pi _G)$ and $G$ acts on $G^*$ by the dressing action in the sense of Lu and Weinstein \cite{LW}. Then $\sigma $ is a quasi-Poisson-Hamiltonian action. Actually, by the definition, the map $\mu $ satisfies
\begin{eqnarray}\label{P-mm}
\pi ^\sharp (\mu ^*(X^L))=X_M
\end{eqnarray}
for any $X$ in $\mathfrak{g}$, where $X^L$ is a left-invariant $1$-form on $G^*$ with value $X$ at $e$ in $G^*$. The quotient manifold $D/G$ is diffeomorphic to $G^*$ as a manifold. The quasi-triple $(D,G,\mathfrak{g}^*)$ is complete since $\pi _G$ is complete. Then $1$-form $\hat{X}_{\mathfrak{g}^*}$ defined by (\ref{1-form}) is global for any $X$ in $\mathfrak{g}$. Furthermore the $1$-form $\hat{X}_{\mathfrak{g}^*}$ on $D/G\cong G^*$ coincides with $X^L$. The complement $\mathfrak{g}^*$ is admissible at any point in $D/G$, so that the map $\mu :M\rightarrow G^*\cong D/G$ is a moment map for the quasi-Poisson action $\sigma $ because of (\ref{P-mm}) and Proposition \ref{sufficient}.
\end{example}

\begin{example}[symplectic manifolds \cite{AK},\cite{Si}]\label{Ex.symplectic}
Let $(M,\omega )$ be a symplectic manifold on which a connected Lie group $G$ acts by a symplectic-Hamiltonian action $\sigma $. Since the symplectic structure $\omega $ induces a Poisson structure $\pi$, the pair $(M,\pi )$ is a Poisson manifold. Then the action $\sigma $ is a Poisson action of a trivial Poisson Lie group $(G,0)$ on $(M,\pi )$. The trivial Poisson structure $0$ on $G$ is complete and the quasi-triple corresponding to $(G,0)$ is $(T^*G,G,\mathfrak{g}^*)$, where $T^*G\cong G\times \mathfrak{g}^*$ is the cotangent bundle of $G$ equipped with the group structure of a semi-direct product with respect to coadjoint action of $G$ on $\mathfrak{g}^*$ (See \cite{AK}). The dual group $G^*$ of $(G,0)$ is the additive group $\mathfrak{g}^*$ and the moment map $\mu $ for symplectic action $\sigma $ is $G$-equivariant with respect to $\sigma $ on $M$ and $\mbox{Ad}^*$ on $\mathfrak{g}^*$ by the definition. Furthermore the dressing action of $G$ on $G^*=\mathfrak{g}^*$ coincides with the coadjoint action $\mathrm{Ad}^*$. Thus the map $\mu :M\rightarrow \mathfrak{g}^*=G^*$ is a moment map for the Poisson action $\sigma $. Therefore, by Example \ref{Ex.Poisson}, the map $\mu :M\rightarrow \mathfrak{g}^*\cong T^*G/G$ is a moment map for the quasi-Poisson action $\sigma $ on the quasi-Poisson $(G,0,0)$-manifold $(M,\pi )$.
\end{example}

\section{Main Result}\label{Main Result}

In this section, we carry out deformations of symplectic structures on a smooth manifold. We use the moment map theory for quasi-Poisson actions for it. A moment map for the quasi-Poisson action on a quasi-Poisson $G_D^\mathfrak{h}$-manifold $(M,\pi _M^\mathfrak{h})$ are defined with the conditions for the family of quasi-Poisson $G_D^{\mathfrak{h}'}$-structures $\left\{\pi _M^{\mathfrak{h}'}\right\}_{\mathfrak{h}\rq{}}$ on $M$. For each complement $\mathfrak{h}\rq{}$, there exists a twist $t$ in $\Lambda ^2\mathfrak{g}$ such that $\mathfrak{h}\rq{}=\mathfrak{h}_t$, so that the family $\left\{\pi _M^{\mathfrak{h}'}\right\}_{\mathfrak{h}\rq{}}$ is regarded as the family parametrized by twist, $\left\{\pi _M^{\mathfrak{h}_t}\right\}_{t\in \Lambda ^2\mathfrak{g}}$. %The conditions for moment maps for quasi-Poisson actions are imposed on the family of quasi-Poisson $G_D^{\mathfrak{h}_t}$-structures parametrized by twists $t$ in $\Lambda ^2\mathfrak{g}$. 
When the quasi-Poisson $G_D^{\mathfrak{h}_t}$-structure with twist $t=0$ is induced by a given symplectic structure, we will give the method of finding a quasi-Poisson $G_D^{\mathfrak{h}_t}$-structure which induced a symplectic structure in $\left\{\pi _M^{\mathfrak{h}_t}\right\}_t$. That is, we can deform a given symplectic structure to a new one by a twist $t$.
%twists $t$ by which the quasi-Poisson $G_D^{\mathfrak{h}_t}$-structures parametrized induce symplectic structures in  $\left\{\pi _M^{\mathfrak{h}_t}\right\}_{t\in \Lambda ^2\mathfrak{g}}$. %we call it a {\it deformation of the symplectic structure by twist $t$} to specify a quasi-Poisson $G_D^{\mathfrak{h}_t}$-structure inducing a symplectic structure or the induced symplectic structure itself. 
This deformation can be carried out due to using the family $\left\{\pi _M^{\mathfrak{h}_t}\right\}_t$ as moment map conditions for quasi-Poisson actions. In this regard, it is described as follows in \cite{AK}: It would be interesting to find a geometric framework for considering the family $\left\{\pi _M^{\mathfrak{h}_t}\right\}_t$. Our deformation is one of the answers for this proposal.
%We consider that this method is one of geometry frameworks for using the family. It is described that considering the framework is interesting, in \cite{AK}.

Let $(M,\omega )$ be a symplectic manifold on which an $n$-dimensional connected Lie group $G$ acts by symplectic-Hamiltonian action $\sigma $ with a moment map $\mu :M\rightarrow \mathfrak{g}^*$. Let $\pi$ be the nondegenerate Poisson structure on $M$ induced by $\omega $. Then $\mu $ is a moment map for the quasi-Poisson-Hamiltonian action $\sigma $ of $(G,0,0)$ on $(M,\pi )$ by Example \ref{Ex.symplectic} in Section \ref{quasi-Poisson}. 

Let $(\mathfrak{g}\oplus \mathfrak{g}^*,\mathfrak{g},\mathfrak{g}^*)$ be the Manin triple corresponding to the trivial Poisson Lie group $(G,0)$, where $\mathfrak{g}\oplus \mathfrak{g}^*$ has the Lie bracket
\begin{eqnarray}
[X,Y]=[X,Y]_\mathfrak{g},\quad [X,\xi ]=\mathrm{ad}^*_X\xi ,\quad [\xi ,\eta ]=[\xi ,\eta ]_{\mathfrak{g}^*}=0
\end{eqnarray}
for any $X, Y$ in $\mathfrak{g}$ and $\xi ,\eta $ in $\mathfrak{g}^*$.
Here the bracket $[\cdot ,\cdot ]_\mathfrak{g}$ and $[\cdot ,\cdot ]_{\mathfrak{g}^*}$ are the brackets on $\mathfrak{g}$ and $\mathfrak{g}^*$ respectively. Then the Manin (quasi-)triple $(\mathfrak{g}\oplus \mathfrak{g}^*,\mathfrak{g},\mathfrak{g}^*)$ defines $F:=F_{\mathfrak{g}^*}=0$ and $\varphi :=\varphi _{\mathfrak{g}^*}=0$ (see (\ref{varphi}) and (\ref{cobracket dual})). Since the corresponding quasi-triple $(T^*G,G,\mathfrak{g}^*)$ is complete by Example \ref{Ex.Poisson} and \ref{Ex.symplectic}, an isotropic complement $\mathfrak{g}^* $ is admissible at any $\xi $ in $ \mathfrak{g}^*$ by Definition \ref{admissible}, and hence it is admissible at any $\xi $ in $ \mu (M)$. 

Let $\mathfrak{g}^*_t$ be an isotropic complement of $\mathfrak{g}$ in $\mathfrak{g}\oplus \mathfrak{g}^*$ with a twist $t$ in $\Lambda ^2\mathfrak{g}$ from $\mathfrak{g}^*$. When we deform $\pi $ to $\pi _M^t:=\pi -t_M$ by a twist $t$, the quasi-Poisson Lie group $(G,0,0)$ is deformed to $(G,\pi _G^t, \varphi_{\mathfrak{g}^*_t})$, where $\pi _G^t=t^L-t^R$ and $\varphi_{\mathfrak{g}^*_t}=\frac{1}{2}[t,t]+\varphi _t$ by (\ref{phi-t}) and (\ref{pi_G-t}). Moreover it follows from $F=0$ and (\ref{components}) that $\varphi _t=0$. So $\varphi_{\mathfrak{g}^*_t}=\frac{1}{2}[t,t]$.

On the other hand, it follows from Definition \ref{q-P-mfd} that the quasi-Poisson $(G,\pi _G^t, \varphi_{\mathfrak{g}^*_t})$-manifold $(M,\pi _M^t)$ satisfies
\begin{eqnarray}
\frac{1}{2}\left[\pi _M^t,\pi _M^t\right]&=&(\varphi_{\mathfrak{g}^*_t})_M, \label{q-P-a-3}\\
\mathfrak{L}_{X_M}\pi _M^t&=&F_{\mathfrak{g}^*_t}(X)_M.\label{q-P-a-4}
\end{eqnarray}
If $(\varphi_{\mathfrak{g}^*_t})_M=0$, i.e., $[t,t]_M=0$, then the $2$-vector field $\pi _M^t$ is a Poisson structure on $M$ by (\ref{q-P-a-3}). 

Assume that a twist $t$ in $\Lambda^2\mathfrak{g}$ is an {\it r-matrix}, namely that $\left[t,t\right]$ is $\mbox{ad}$-invariant. Then $\pi _G^t=t^L-t^R$ is a multiplicative Poisson structure (see \cite{LW}). Therefore $(G,\pi _G^t)$ is a Poisson Lie group. Then it follows that $F_{\mathfrak{g}^*_t}$ coincides with the dual of the bracket map $[\cdot ,\cdot ]_{\pi _G^t}\!:\mathfrak{g}^*\wedge \mathfrak{g}^*\rightarrow \mathfrak{g}^*$ on $\mathfrak{g}^*$ defined by the Poisson Lie group $(G,\pi _G^t)$. In fact, by the relation (\ref{P-mm}), we have
\begin{eqnarray}\label{coincide with the bracket}
F_{\mathfrak{g}^*_t}^*(\xi ,\eta )=\mathrm{ad}_{t^\sharp (\xi )}\eta -\mathrm{ad}_{t^\sharp (\eta )}\xi ,
\end{eqnarray}
where $<t^\sharp (\xi ),\eta >:=t(\xi ,\eta )$. And the bracket on $\mathfrak{g}^*$ induced by a multiplicative Poisson structure defined by an r-matrix is represented by the right-hand side of (\ref{coincide with the bracket}) (see \cite{L2}, Ex.2.19). Therefore, since $G$ is connected, the condition (\ref{q-P-a-4}) means that the action $\sigma $ is a Poisson action of $(G,\pi _G^t)$ on $(M,\pi _M^t)$ under the assumption that $t$ is an r-matrix and that $[t,t]_M=0$. 

Let $\{e_i\}$ be a basis on $\mathfrak{g}$, a set $\{\varepsilon ^i\}$ be the dual basis of $\{e_i\}$ on $\mathfrak{g}^*$. Then we can write by (\ref{basis}),
\begin{eqnarray}
\mathfrak{g}^*_t=\mathrm{span}\{\varepsilon ^i+t^{ij}e_j\ |\ i=1,\dots ,n\},
\end{eqnarray}
where $t=\frac{1}{2}t^{ij}e_i\wedge e_j$ in $\Lambda ^2\mathfrak{g}$. If $\mathfrak{g}^*_t$ is admissible at any point in $\mu (M)$, then it satisfies $\mathrm{Im}\pi _p^\sharp =\mathrm{Im}(\pi _M^t)_p^\sharp $ for any $p$ in $ M$ by Theorem \ref{generalized foliation}. The nondegeneracy of $\pi$ means that $\mathrm{Im}\pi _p^\sharp =T_pM$ for any $p$ in $M$. Therefore, by the fact that $\mathrm{Im}(\pi _M^t)_p^\sharp =T_pM$ for any $p$ in $M$, a quasi-Poisson structure $\pi _M^t$ is also nondegenerate.

Here we shall examine the condition for a isotropic complement to be admissible at a point in $\mathfrak{g}^*$ in more detail. Let $(\xi _i)$ be the linear coordinates for $\{\varepsilon^i\}$. Then it follows that for $i=1,\dots ,n$,
\begin{eqnarray}\label{infinitesimal action}
(\varepsilon ^i+t^{ij}e_j)_{\mathfrak{g}^*}&=&-\frac{\partial}{\partial \xi _i}+t^{ij}c_{jl}^k\xi_k\frac{\partial}{\partial \xi _l}\nonumber\\
                            &=&-t^{ij}\sum_{l\neq i}^{}c_{lj}^k\xi_k\frac{\partial}{\partial \xi _l}-(1+t^{ij}c_{ij}^k\xi_k)\frac{\partial}{\partial \xi _i},
\end{eqnarray}
where $X\mapsto X_{\mathfrak{g}^*}$, for $X$ in $\mathfrak{g}\oplus \mathfrak{g}^*$, is the infinitesimal action of the dressing action on $\mathfrak{g}^*\cong T^*G/G$. The quasi-triple $(T^*G,G,\mathfrak{g}_t^*)$ is complete if and only if the elements (\ref{infinitesimal action}) form a basis on $T_\xi (\mathfrak{g}^*)\cong \mathfrak{g}^*$ for any $\xi=(\xi_1,\dots ,\xi_n)$. Hence this means that the matrix
\begin{eqnarray}\label{admissible matrix}
A_t(\xi ):=\left( 
\begin{matrix}
  -1-t^{1j}c_{1j}^k\xi_k & -t^{1j}c_{2j}^k\xi_k & \cdots  & -t^{1j}c_{nj}^k\xi_k\\
  -t^{2j}c_{1j}^k\xi_k & -1-t^{2j}c_{2j}^k\xi_k & \cdots  & -t^{2j}c_{nj}^k\xi_k\\
  \vdots & \vdots & \ddots & \vdots \\
  -t^{nj}c_{1j}^k\xi_k & -t^{nj}c_{2j}^k\xi_k & \cdots  & -1-t^{nj}c_{nj}^k\xi_k
\end{matrix}
\right)
\end{eqnarray}
is regular for any $\xi $. Therefore this is equivalent to $f_t(\xi _1,\dots ,\xi _n):=\det A_t(\xi )\neq 0$. Since the constant term of $f_t$ is $(-1)^n$ and since coefficients of the rest of the term include $t^{ij}$'s, a family of hypersurfaces $\{f_t=0\}_{t\in \Lambda ^2\mathfrak{g}}$ in $\mathfrak{g}^*$ diverges to infinity as $t$ approaches the origin $0$ in $\Lambda ^2\mathfrak{g}$. If $M$ is compact, then $\mu (M)$ is bounded. So it follows that an intersection of $\{f_t=0\}$ and $\mu (M)$ is empty for a twist $t$ close sufficiently to the origin. Therefore since $\mathfrak{g}^*_t$ is admissible on $\mu (M)$, the $2$-vector field $\pi _M^t$ is nondegenerate.

Since any nondegenerate Poisson structure on $M$ defines a symplectic structure on $M$, the following theorem holds. 

\begin{theorem}\label{abstract main}
Let $(M,\omega)$ be a symplectic manifold on which a connected Lie group $G$ with the Lie algebra $\mathfrak{g}^*$ acts by a symplectic-Hamiltonian action $\sigma $ and $\mu $ be a moment map for $\sigma $. Then the following holds:
\begin{enumerate}
\item If a twist $t$ in $\Lambda ^2\mathfrak{g}$ satisfies that $[t,t]_M=0$, then $t$ deforms the Poisson structure $\pi $ induced by $\omega $ to a Poisson structure $\pi _M^t:=\pi -t_M$. Moreover, if $t$ is an r-matrix, then $\sigma $ is a Poisson action of $(G,\pi _G^t)$ on $(M,\pi _M^t)$, where $\pi _G^t=t^L-t^R$. \label{(i)}
\item For a twist $t$ in $\Lambda ^2\mathfrak{g}$, if the isotropic complement $\mathfrak{g}_t^*$ is admissible on $\mu (M)$, then $t$ deforms the nondegenerate $2$-vector field $\pi $ induced by $\omega $ to a nondegenerate $2$-vector field $\pi _M^t$. In particular, if $M$ is compact, then a $2$-vector field $\pi _M^t$ is nondegenerate for a twist $t$ close sufficiently to the origin $0$ in $\Lambda ^2\mathfrak{g}$. \label
{(ii)}
\end{enumerate}
Therefore, if a twist $t$ satisfies the assumptions of both \ref{(i)} and \ref{(ii)}, then $t$ deforms $\omega $ to a symplectic structure $\omega^t$ induced by the nondegenerate Poisson structure $\pi _M^t$.
\end{theorem}

The following theorem gives a sufficient condition for a twist to deform a symplectic structure.

\begin{theorem}\label{main}
Let $(M,\omega)$ be a symplectic manifold on which an $n$-dimensional connected Lie group $G$ acts by a symplectic-Hamiltonian action $\sigma $. Assume that $X,Y$ in $\mathfrak{g}$ satisfy $[X,Y]=0$. Then the twist $t=\frac{1}{2}X\wedge Y$ in $\Lambda ^2\mathfrak{g}$ deforms the symplectic structure $\omega $ to a symplectic structure $\omega_t$. For example, a twist $t$ in $ \Lambda ^2\mathfrak{h}$, where $\mathfrak{h}$ is a Cartan subalgebra of $\mathfrak{g}$, satisfies the assumption of the theorem.
\end{theorem}

\begin{proof}
For $X$ and $Y$ in $ \mathfrak{g}$, we set
$$X=X^ie_i,\ Y=Y^je_j,$$
where $\{e_i\}_{i=1}^n$ is a basis on the Lie algebra $\mathfrak{g}$. Then since $[X,Y]=X^iY^jc_{ij}^ke_k=0$, we obtain the following conditions:
\begin{eqnarray*}
X^iY^jc_{ij}^k=0\quad \text{for any $k$},
\end{eqnarray*}
where $c_{ij}^k$ are the structure constants of $\mathfrak{g}$ with respect to the basis $\{e_i\}$. Moreover, since we have
\begin{eqnarray*}
[t,t]=\left[\frac{1}{2}X\wedge Y,\frac{1}{2}X\wedge Y\right]=\frac{1}{2}X\wedge [X,Y]\wedge Y=0,
\end{eqnarray*}
the twist $t$ is an r-matrix such that $[t,t]_M=0$ obviously. Hence $\pi _M^t$ is a Poisson structure, and if $\pi _M^t$ is nondegenerate, then twist $t$ induces the symplectic structure $\omega _t$.

We shall show the nondegeneracy of $\pi _M^t$. Let $\mu $ be the moment map for a given symplectic-Hamiltonian action $\psi $. Then the nondegeneracy of $\pi _M^t$ means that $\mathfrak{g}^*_t$ is admissible at any point in $\mu (M)$. We prove a stronger condition that the quasi-triple $(T^*G,G,\mathfrak{g}^*_t)$ is complete.

Let $\{\varepsilon ^i\}$ be the dual basis of $\{e_i\}$ on $\mathfrak{g}^*$ and $(\xi _i)$ be the linear coordinates for $\{\varepsilon ^i\}$. Since $t=\frac{1}{2}X^iY^je_i\wedge e_j$,
$$\mathfrak{g}^*_t=\mbox{span}\{\varepsilon ^i+X^iY^je_j|i=1,\cdots n\}.$$
Then it follows that for $i=1,\dots ,n$,
\begin{eqnarray}\label{infinitesimal action-2}
(\varepsilon ^i+X^iY^je_j)_{\mathfrak{g}^*}%&=&-\frac{\partial}{\partial \xi _i}+X^iY^jc_{jl}^k\xi_k\frac{\partial}{\partial \xi _l}\nonumber\\
                            &=&-X^iY^j\sum_{l\neq i}^{}c_{lj}^k\xi_k\frac{\partial}{\partial \xi _l}-(1+X^iY^jc_{ij}^k\xi_k)\frac{\partial}{\partial \xi _i}.
\end{eqnarray}
The quasi-triple $(T^*G,G,\mathfrak{g}_t^*)$ is complete if and only if the elements (\ref{infinitesimal action-2}) form a basis on $T_\xi (\mathfrak{g}^*)\cong \mathfrak{g}^*$ for any $\xi=(\xi_1,\dots ,\xi_n)$. Therefore we shall prove that the matrix
\begin{eqnarray}\label{matrix-2}
\left( 
\begin{matrix}
  -1-X^1Y^jc_{1j}^k\xi_k & -X^1Y^jc_{2j}^k\xi_k & \cdots  & -X^1Y^jc_{nj}^k\xi_k\\
  -X^2Y^jc_{1j}^k\xi_k & -1-X^2Y^jc_{2j}^k\xi_k & \cdots  & -X^2Y^jc_{nj}^k\xi_k\\
  \vdots & \vdots & \ddots & \vdots \\
  -X^nY^jc_{1j}^k\xi_k & -X^nY^jc_{2j}^k\xi_k & \cdots  & -1-X^nY^jc_{nj}^k\xi_k
\end{matrix}
\right)
\end{eqnarray}
is regular. In the case of $X=0$, this matrix is equal to the opposite of the identity matrix, so that it is regular. In the case of $X\neq 0$,
using $X^iY^jc_{ij}^k=0$, we can transform the matrix to the opposite of the identity matrix.
Thus the matrix (\ref{matrix-2}) is regular. Therefore $\mathfrak{g}^*_t$ is admissible at any point in $\mathfrak{g}^*$. That is, $(T^*G,G,\mathfrak{g}^*_t)$ is complete.
\end{proof}

\begin{rem}
We try to generalize the assumption of Theorem \ref{main} and consider $X,Y$ in $ \mathfrak{g}$ such that $[X,Y]=aX+bY\ (a,b$ in $ \mathbb{R})$, that is, the subspace spanned by $X,Y$ is also a Lie subalgebra. We set $t=\frac{1}{2}X\wedge Y$ in $\Lambda ^2\mathfrak{g}$. Since
\begin{eqnarray*}
[t,t]=\frac{1}{2}X\wedge [X,Y]\wedge Y=\frac{1}{2}X\wedge (aX+bY)\wedge Y=0,
\end{eqnarray*}
the twist $t$ is an r-matrix such that $[t,t]_M=0
$. Therefore the symplectic action $\psi $ is a Poisson action of $(G,\pi _G^t)$ on $(M, \pi _M^t)$. Then we research whether $\mathfrak{g}_t^*$ is admissible at all points in $\mathfrak{g}^*$. Similarly to the proof of Theorem \ref{main}, a matrix to check the regularity can be deformed to
$$\left( 
\begin{matrix}
  -1-(aX^k+bY^k)\xi _k & 0 & \cdots  & 0\\
  0 & -1 & \cdots  & 0\\
  \vdots & \vdots & \ddots & \vdots \\
  0 & 0 & \cdots  & -1
\end{matrix}
\right).$$
Therefore this matrix is regular if and only if 
$$-1-(aX^k+bY^k)\xi _k\neq 0.$$
In the case of $[X,Y]=0$, by Theorem \ref{main}, the space $\mathfrak{g}_t^*$ is admissible at all points in $\mathfrak{g}^*$. In the case of $[X,Y]\neq 0$, we shall denote by $\cdot $ the standard inner product on $\mathfrak{g}^*\cong \mathbb{R}^n$. Then the above condition means
$$\eta _{[X,Y]}\cdot \xi \neq -1,$$
where $\eta _X=\sum_k^{}X^k\varepsilon ^k$ for $X=X^ke_k$ in $\mathfrak{g}$. Let $\xi '$ be an element which is not orthogonal to $\eta _{[X,Y]}$. By setting
$$\xi :=-\frac{\xi '}{\eta _{[X,Y]}\cdot \xi '},$$
we obtain $\eta _{[X,Y]}\cdot \xi =-1$, so that $\mathfrak{g}_t^*$ is not admissible at $\xi$. Eventually, to make sure of the admissibility of $\mathfrak{g}_t^*$, we need check whether such a point $\xi $ is included in $\mu (M)$.
\end{rem}

\section{Deformations of the canonical symplectic form on $\mathbb{C}\mathrm{P}^n$ and $\mathrm{Gr}(r;\mathbb{C}^n)$}

 In this section, we compute specifically which element $t$ in $ \Lambda ^2\mathfrak{g}$ defines a different symplectic structure $\omega _t$ from given one $\omega $ on a smooth manifold.
 One example is the complex projective space $(\mathbb{C}\mathrm{P}^k,\omega _{\mathrm{FS}})$, where $\omega _{\mathrm{FS}}$ is the Fubini-Study form, with an action of $\mathrm{SU}(k+1)$, $k=1,2$. Another is $(\mathbb{C}\mathrm{P}^n,\omega _{\mathrm{FS}})$ with an action of the torus group $\mathbb{T}^n$. The other is the complex Grassmannian $(\mathrm{Gr}(r;\mathbb{C}^n),\omega _{\mathrm{KK}})$, where $\omega _{\mathrm{KK}}$ is the Kirillov-Kostant form, with an action of $\mathrm{SU}(k+1)$.

 First we review the relation between $\mathrm{SU}(n+1)$ and $\mathbb{C}\mathrm{P}^n$. For any $\left[z_1:\dots :z_{n+1}\right]$ in $\mathbb{C}\mathrm{P}^n$ and $g=(a_{ij})$ in $\mathrm{SU}(n+1)$, the action is given by
$$g\cdot \left[z_1:\dots :z_{n+1}\right]:=\left[\sum_{j=1}^{n+1}a_{1j}z_j: \dots :\sum_{j=1}^{n+1}a_{n+1,j}z_j\right].$$
The isotropic subgroup of $\left[1:0:\dots :0\right]$ is
$$\mathrm{S}(\mathrm{U}(1)\times \mathrm{U}(n))=\left\{ \left.
\left( 
\begin{matrix}
 e^{i\theta} & O \\
 O           & B 
\end{matrix}
\right)
\in  \mathrm{SU}(n+1)\right|\theta \in \mathbb{R},B\in U(n)\right\}.$$
Therefore it follows
$$\mathrm{SU}(n+1)/\mathrm{S}(\mathrm{U}(1)\times \mathrm{U}(n))\cong \mathbb{C}\mathrm{P}^n.$$
The complex projective space $\mathbb{C}\mathrm{P}^n$ has the coordinate neighborhood system $\{(U_i,\varphi _i)\}_i$ consisting of $n+1$ open sets $U_i$ given by
\begin{eqnarray*}
&U_i&:=\{\left[z_1:\dots :z_{n+1}\right]\in \mathbb{C}\mathrm{P}^n|z_i\neq 0\},\\
&\varphi _i&:U_i\rightarrow \mathbb{C}^n\cong \mathbb{R}^{2n},\\
\end{eqnarray*}
\vspace{-13mm}
\begin{eqnarray*}
\left[z_1:\dots :z_{n+1}\right]&\mapsto &\left(\frac{z_1}{z_i},\dots ,\frac{z_{i-1}}{z_i},\frac{z_{i+1}}{z_i},\dots ,\frac{z_{n+1}}{z_i}\right)\\
&\mapsto &\left(\mbox{Re}\frac{z_1}{z_i},\mbox{Im}\frac{z_1}{z_i},\dots ,\mbox{Re}\frac{z_{n+1}}{z_i},\mbox{Im}\frac{z_{n+1}}{z_i}\right),
\end{eqnarray*}
for $i=1,\dots ,n+1$. By using this coordinate system, the Fubini-Study form $\omega _{\mathrm{FS}}$ on $\mathbb{C}\mathrm{P}^n$ is defined by setting
$$\varphi _i^*\left(\frac{i}{2}\partial \bar{\partial }\log \left(\sum_j^{}|z_j|^2+1\right)\right)$$
on each $U_i$.

 The action of $\mathrm{SU}(n+1)$ on $(\mathbb{C}\mathrm{P}^n,\omega _{\mathrm{FS}})$ is a symplectic-Hamiltonian action and its moment map $\mu $ satisfies
\begin{eqnarray*}
&<\mu (\left[z_1:\dots :z_{n+1}\right]),X>&=\frac{1}{2}\mbox{Im}\frac{\langle {}^t \!(z_1,\dots ,z_{n+1}),X{}^t \!(z_1,\dots ,z_{n+1})\rangle}{\langle{}^t \!(z_1,\dots ,z_{n+1}),{}^t \!(z_1,\dots ,z_{n+1})\rangle}\\
\end{eqnarray*}
for any $\left[z_1:\dots :z_{n+1}\right]$ in $\mathbb{C}\mathrm{P}^n$ and $X$ in $\mathfrak{su}(n+1)$. 
%Remark that $(\mathbb{C}\mathrm{P}^n,\omega _{\mathrm{FS}})$ is identified with $(\mathcal{O},\omega _{\mathcal{O}})$, where $\mathcal{O}$ is a coadjoint orbit of $\mathrm{SU}(n+1)$ on $\mathfrak{su}(n)^*$ and $\omega _{\mathcal{O}}$ is the Kostant-Kirillov form on $\mathcal{O}$. Then $\mu $ is identified with the moment map for the coadjoint action
%$$\mu _{\mathcal{O}}:\mathcal{O}\hookrightarrow \mathfrak{su}(n+1)^*,\ \xi \mapsto \xi $$
%(see e.g., ).

We use
\begin{eqnarray*}
&X_{ij}&:\text{the $(i,j)$-element is $1$, the $(j,i)$-element is $-1$, and the rest are $0$},\\
&Y_{ij}&:\text{the $(i,j)$- and $(j,i)$-elements are $i$, and the rest are $0$},\\
&Z_k&:\text{the $(k,k)$--element is $i$, and the $(n+1,n+1)$-element is $-i$}
\end{eqnarray*}
for $1\leq i<j\leq n+1$ and $k=1,\dots , n$, as a basis of $\mathfrak{su}(n+1)$ which is defined by a Chevalley basis of the complexified Lie algebra $\mathfrak{sl}(n+1,\mathbb{C})$ of $\mathfrak{su}(n+1)$. The subspace spanned by $Z_k$'s is a Cartan subalgebra of $\mathfrak{su}(n+1)$.

We consider the case of $n=1$. The complex projective line $\mathbb{C}\mathrm{P}^1$ has the coordinate neighborhood system $\{(U_1,\varphi _1),(U_2,\varphi _2)\}$ given by
\begin{eqnarray*}
&U_i&:=\{\left[z_1:z_2\right]\in \mathbb{C}\mathrm{P}^1|z_i\neq 0\}\ (i=1,2),\\
&\varphi _1&:U_1\rightarrow \mathbb{C}\cong \mathbb{R}^2,\ \left[z_1:z_2\right]\mapsto \frac{z_2}{z_1}\mapsto \left(\mbox{Re}\frac{z_2}{z_1},\mbox{Im}\frac{z_2}{z_1}\right),\\
&\varphi _2&:U_2\rightarrow \mathbb{C}\cong \mathbb{R}^2,\ \left[z_1:z_2\right]\mapsto \frac{z_1}{z_2}\mapsto \left(\mbox{Re}\frac{z_1}{z_2},\mbox{Im}\frac{z_1}{z_2}\right).
\end{eqnarray*}
The Fubini-Study form $\omega _{\mathrm{FS}}$ on $\mathbb{C}\mathrm{P}^1$ is
$$\omega _{\mathrm{FS}}=\frac{dx_1\wedge dy_1}{(x_1^2+y_1^2+1)^2}$$
on $U_1$, where $(x_1,y_1):=\left(\mathrm{Re}\frac{z_2}{z_1},\mathrm{Im}\frac{z_2}{z_1}\right)$. Then a moment map $\mu :\mathbb{C}\mathrm{P}^1\rightarrow \mathfrak{su}(2)^*$ for the natural action of $\mathrm{SU}(2)$ on $(\mathbb{C}\mathrm{P}^1,\omega _{\mathrm{FS}})$ is defined by 
$$<\mu (\left[z_1:z_2\right]),X>=-\frac{1}{2}\mbox{Im}\frac{\langle {}^t \!(z_1,z_2),X{}^t \!(z_1,z_2)\rangle}{\langle{}^t \!(z_1,z_2),{}^t \!(z_1,z_2)\rangle}$$
for any $\left[z_1:z_2\right]$ in $\mathbb{C}\mathrm{P}^1$ and $X$ in $\mathfrak{su}(2)$. Then $e_1:=X_{12}$, $e_2:=Y_{12}$ and $e_3:=Z_1$
%%%%%%%%%%%%%%%%%%
\iffalse
$$e_1:=X_{12}=
\left( 
\begin{matrix}
 0  & 1 \\
 -1 & 0 
\end{matrix}
\right),\ 
e_2:=Y_{12}=
\left( 
\begin{matrix}
 0 & i \\
 i & 0 
\end{matrix}
\right),\ 
e_3:=Z_1=
\left( 
\begin{matrix}
 i & 0 \\
 0 & -i 
\end{matrix}
\right)$$
\fi
%%%%%%%%%%%%%%%%%
form a basis of $\mathfrak{su}(2)$. Let $\{\varepsilon ^i\}$ be the dual basis of $\mathfrak{su}(2)^*$. We obtain
\begin{eqnarray*}
\mu (x_1,y_1)&=&\frac{y_1}{1+x_1^2+y_1^2}\varepsilon^1+\frac{x_1}{1+x_1^2+y_1^2}\varepsilon^2+\frac{1-x_1^2-y_1^2}{2(1+x_1^2+y_1^2)}\varepsilon^3.
\end{eqnarray*}
Hence $\mu (\mathbb{C}\mathrm{P}^1)\subset \mathfrak{su}(2)^*$ is the $2$-sphere with center at the origin and with radius $\frac{1}{2}$.

 Let $(\xi _i)$ be the linear coordinates for $\{\varepsilon ^i\}$. We set $\mathfrak{g}:=\mathfrak{su}(2)$. Any twist $t$ is an r-matrix on $\mathfrak{g}$ because $e_1\wedge e_2\wedge e_3$ is $\mathrm{ad}$-invariant. Since $\mathbb{C}\mathrm{P}^1$ is $2$-dimensional, it follows that $[t,t]_{\mathbb{C}\mathrm{P}^1}=0$. Therefore we can deform the Poisson structure $\pi _{\mathrm{FS}}$ induced by $\omega _{\mathrm{FS}}$ to a Poisson structure $\pi _{\mathrm{FS}}^t$ on $\mathbb{C}\mathrm{P}^1$ by $t$ 
% $2$-vector field $\pi _{\mathrm{FS}}^t$ deformed $\pi _{\mathrm{FS}}$, which is the Poisson structure on $\mathbb{C}\mathrm{P}^1$ defined by $\omega _{\mathrm{FS}}$, by $t$ is a Poisson structure on $\mathbb{C}\mathrm{P}^1$ 
and the natural action is a Poisson action of $(\mathrm{SU}(2),t^L-t^R)$.

Let $\mathfrak{g}^*_t$ be the space twisted $\mathfrak{g}^*$ by $t$ in $ \Lambda ^2\mathfrak{g}$. We consider what is the condition for $t$ under which $\mathfrak{g}^*_t$ is admissible on $\mu (\mathbb{C}\mathrm{P}^1)$. For any twist
$$t=\sum _{i<j}^{}\frac{1}{2}\lambda _{ij}e_i\wedge e_j\in \Lambda ^2\mathfrak{g}\ (\lambda _{ij}\in \mathbb{R}),$$
we obtain
$$\mathfrak{g}^*_t=\mbox{span}\{\varepsilon ^1+\lambda _{12}e_2+\lambda _{13}e_3,\ \varepsilon ^2-\lambda _{12}e_1+\lambda _{13}e_3,\ \varepsilon ^3-\lambda _{13}e_1-\lambda _{23}e_2\}.$$
We calculate as 
\begin{eqnarray*}
(\varepsilon ^1+\lambda _{12}e_2+\lambda _{13}e_3)_{\mathfrak{g}^*}&=&-(1+2\lambda _{12}\xi _3-2\lambda _{13}\xi _2)\frac{\partial }{\partial \xi_1}-2\lambda _{13}\xi _1\frac{\partial }{\partial \xi_2}+2\lambda _{12}\xi _1\frac{\partial }{\partial \xi_3},\\
(\varepsilon ^2-\lambda _{12}e_1+\lambda _{13}e_3)_{\mathfrak{g}^*}&=&2\lambda _{23}\xi _2\frac{\partial }{\partial \xi_1}-(1+2\lambda _{12}\xi _3+2\lambda _{23}\xi _1)\frac{\partial }{\partial \xi_2}+2\lambda _{12}\xi _2\frac{\partial }{\partial \xi_3},\\
(\varepsilon ^3-\lambda _{13}e_1-\lambda _{23}e_2)_{\mathfrak{g}^*}&=&2\lambda _{23}\xi _3\frac{\partial }{\partial \xi_1}-2\lambda _{13}\xi _3\frac{\partial }{\partial \xi_2}-(1-2\lambda _{13}\xi _2+2\lambda _{23}\xi _1)\frac{\partial }{\partial \xi_3}.
\end{eqnarray*}
Then $\mathfrak{g}^*_t$ is admissible at $\xi =(\xi_1,\xi_2,\xi_3)$ in $\mathfrak{g}^*$ if and only if %$T_\xi (\mathfrak{g}^*)$ is spanned by
%$$\{(\varepsilon ^1+\lambda _{12}e_2+\lambda _{13}e_3)_{\mathfrak{g}^*,\xi },(\varepsilon ^2-\lambda _{12}e_1+\lambda _{13}e_3)_{\mathfrak{g}^*,\xi },(\varepsilon ^3-\lambda _{13}e_1-\lambda _{23}e_2)_{\mathfrak{g}^*,\xi }\}.$$
%Moreover this condition is equivalent to that
the matrix
$$A_t(\xi )=\left( 
\begin{matrix}
 1+2\lambda _{12}\xi _3-2\lambda _{13}\xi _2 & 2\lambda _{13}\xi _1                      & -2\lambda _{12}\xi _1                        \\
 -2\lambda _{23}\xi _2                       & 1+2\lambda _{12}\xi _3+2\lambda _{23}\xi _1 & -2\lambda _{12}\xi _2                       \\
 -2\lambda _{23}\xi _3                       & 2\lambda _{13}\xi _3                      & 1-2\lambda _{13}\xi _2+2\lambda _{23}\xi _1
\end{matrix}
\right)$$
is regular. By computing the determinant of the matrix, we have
\begin{eqnarray*}
f_t(\xi )=\det A_t(\xi )=(1+2\lambda _{23}\xi _1-2\lambda _{13}\xi _2+2\lambda _{12}\xi _3)^2.
\end{eqnarray*}
So the complement $\mathfrak{g}^*_t$ is admissible at $\xi =(\xi_1,\xi_2,\xi_3)$ if and only if $1+2\lambda _{23}\xi _1-2\lambda _{13}\xi _2+2\lambda _{12}\xi _3\neq 0$.

Therefore $\mathfrak{g}^*_t$ is admissible on $\mu (\mathbb{C}\mathrm{P}^1)$ if and only if the "non-admissible surface" $\{\xi =(\xi_1,\xi_2,\xi_3)\in \mathfrak{g}^*|\ 1+2\lambda _{23}\xi _1-2\lambda _{13}\xi _2+2\lambda _{12}\xi _3\neq 0\}$ for $\mathfrak{g}^*_t$ and the image $\mu (\mathbb{C}\mathrm{P}^1)$ have no common point. Since $\mu (\mathbb{C}\mathrm{P}^1)$ is the $2$-sphere with center at the origin and with radius $\frac{1}{2}$, we can see that this condition is equivalent to the condition
$$\lambda _{12}^2+\lambda _{13}^2+\lambda _{23}^2<1.$$
%%%%%%%%%%%%%%%%%%%%%%%%%%%%%%%%%%%%%%%%%%%%%%%%%%%%%%%%%%%%%%%%%%%%%%%%%%%%%%%%%%%%%%%%%%%%%%%%%%%%%%%%%%%%%%%%%%%%%%%%%%%%%%%%%%%%%%%%%%%%%%%%%%%
\iffalse
 for $\vec{t}:={}^t\!(\lambda _{12},\lambda _{13},\lambda _{23})$ to make the solutions of the simultaneous equations\\
\begin{eqnarray}\label{renritsu}
\begin{cases}
1+2\lambda _{23}x-2\lambda _{13}y+2\lambda _{12}z=0 \\
x^2+y^2+z^2=\frac{1}{4}
\end{cases}
\end{eqnarray}
the empty set. By an easy computation, the simultaneous equations (\ref{renritsu}) have no solution if and only if $\parallel \!\vec{t}\!\parallel <1$, where we identify $\Lambda ^2\mathfrak{g}$ with the Euclid space $\mathbb
{R}^3$. %This means that $\vec{t}={}^t\!(\lambda _{12},\lambda _{13},\lambda _{23})$ is included in the open ball with center at the origin and radius $1$ in $\mathbb{R}^3$. 
\fi
%%%%%%%%%%%%%%%%%%%%%%%%%%%%%%%%%%%%%%%%%%%%%%%%%%%%%%%%%%%%%%%%%%%%%%%%%%%%%%%%%%%%%%%%%%%%%%%%%%%%%%%%%%%%%%%%%%%%%%%%%%%%%%%%%%%%%%%%%%%%%%%%%%%
From the above disscusion, we obtain the following theorem.

\begin{theorem}\label{CP1thm}
If a twist $t:=\sum _{i<j}^{}\frac{1}{2}\lambda _{ij}e_i\wedge e_j$ satisfies $\lambda _{12}^2+\lambda _{13}^2+\lambda _{23}^2<1$, then the Fubini-Study form $\omega _{\mathrm{FS}}$ on $\mathbb{C}\mathrm{P}^1$ can be deformed by $t$ in the sense of Section \ref{Main Result}.   
\end{theorem}

We shall see an example of twists giving symplectomorphisms on $\mathbb{C}\mathrm{P}^1$.

\begin{example}\label{cp1}
We use a twist $t=\frac{1}{2}X_{12}\wedge Y_{12}$ in $\Lambda ^2\mathfrak{su}(2)$ and a real number $\lambda $, where $-1<\lambda <1$. The symplectic structure $\omega _{\mathrm{FS}}^{\lambda t}$ deformed $\omega _{\mathrm{FS}}$ by $\lambda t$ is written by
$$\omega _{\mathrm{FS}}^{\lambda t}=\left\{\left(1+\frac{1}{2}\lambda \right)(x_1^2+y_1^2)^2+2(x_1^2+y_1^2)+\left(1-\frac{1}{2}\lambda \right)\right\}^{-1}dx_1\wedge dy_1$$
on $U_1$. Then it follows from an elementary calculation that the symplectic volume $\mathrm{Vol}(\mathbb{C}\mathrm{P}^1, \omega _{\mathrm{FS}}^{\lambda t})$ of $(\mathbb{C}\mathrm{P}^1, \omega _{\mathrm{FS}}^{\lambda t})$ is\\
\begin{eqnarray}\label{taiseki}
\mathrm{Vol}(\mathbb{C}\mathrm{P}^1, \omega _{\mathrm{FS}}^{\lambda t})=
\begin{cases}
\pi &(\lambda =0)\\
\frac{\pi }{\lambda }\log \left|\frac{2+\lambda }{2-\lambda }\right| &(\lambda \neq 0).
\end{cases}
\end{eqnarray}
Next, we consider a cohomology class of each $\omega _{\mathrm{FS}}^{\lambda t}$. Since $H_{\mathrm{DR}}^2(\mathbb{C}\mathrm{P}^1)=\mathbb{R}$, there exists a real number $k_\lambda \in \mathbb{R}$ such that $\left[\omega _{\mathrm{FS}}^{\lambda t}\right]=k_\lambda \left[\omega _{\mathrm{FS}}\right]$. By integrating, we obtain
$$k_\lambda =\frac{1}{\lambda}\log \left|\frac{2+\lambda }{2-\lambda }\right|,$$
where $\lambda \neq 0$. Since the function $k_\lambda $ of $\lambda $ is smooth, %even and strictly monotone increasing when $\lambda $ is positive. Therefore $\omega _{\mathrm{FS}}^{\lambda t}$ and $\omega _{\mathrm{FS}}^{-\lambda t}$ are cohomologous, and 
two symplectic structures $\omega _{\mathrm{FS}}$ and $k_\lambda \omega _{\mathrm{FS}}^{\lambda t}$ are cohomologous. In particular, $(\mathbb{C}\mathrm{P}^1, \omega _{\mathrm{FS}})$ and $(\mathbb{C}\mathrm{P}^1, k_\lambda \omega _{\mathrm{FS}}^{\lambda t})$ are symplectomorphic by Moser's theorem.
\end{example}

Next we shall see deformations of symplectic structures in the case of $\mathbb{C}\mathrm{P}^n$.

\begin{example}
We consider the case of $n=2$. Since $[Y_{23},2Z_1-Z_2]=0$ in $\Lambda ^2\mathfrak{su}(3)$, we use the twist $t=\frac{1}{2}\lambda Y_{23}\wedge (2Z_1-Z_2)\ (\lambda \in \mathbb{R})$ to deform $\omega _{\mathrm{FS}}$. Then $\omega _{\mathrm{FS}}$ is deformed to
\begin{eqnarray*}
\omega_{FS}^t&=&\omega_{FS}+\frac{\lambda}{\left\{\sum_{k}^{}(x_k^2+y_k^2)+1\right\}^3}\left\{(x_1y_2-x_2y_1)dx_1\wedge dy_1\right. \\
&\quad &\quad \quad \quad +(x_1^2-x_2^2)dx_1\wedge dx_2+(x_1y_1-x_2y_2)dx_1\wedge dy_2\\
&\quad &\quad \quad \quad \quad \ +(x_1y_1-x_2y_2)dy_1\wedge dx_2+(y_1^2-y_2^2)dy_1\wedge dy_2\\
&\quad &\quad \quad \quad \quad \quad \quad \quad \quad \quad \quad \quad \quad \quad \quad \ \left. -(x_1y_2-x_2y_1)dx_2\wedge dy_2\right\}
\end{eqnarray*}
on $U_1$, where $x_i:=\mathrm{Re}\frac{z_{i+1}}{z_1}$ and $y_i:=\mathrm{Im}\frac{z_{i+1}}{z_1}$.
\end{example}

\begin{example}
The next example is a symplectic toric manifold $\mathbb{C}\mathrm{P}^n$ with the torus action: 
$$(e^{i\theta _2},e^{i\theta _3},\dots ,e^{i\theta _{n+1}})\cdot [z_1:\dots :z_{n+1}]:=[z_1:e^{i\theta _2}z_2:\dots :e^{i\theta _{n+1}}z_{n+1}]$$
for any $\theta _i$ in $\mathbb{R}$. The moment map $\mu :\mathbb{C}\mathrm{P}^n\rightarrow \mathbb{R}^n$ for this action on $(\mathbb{C}\mathrm{P}^n,\omega _{\mathrm{FS}})$ is
$$\mu ([z_1:\dots :z_{n+1}]):=-\frac{1}{2}\left(\frac{|z_2|^2}{|z|^2},\dots ,\frac{|z_{n+1}|^2}{|z|^2}\right),$$
where $z=(z_1,\dots ,z_{n+1})$ in $\mathbb{C}^n$. We set $X_1:=(1,0,\dots ,0), \dots , X_n:=(0, \dots ,0,1)$. Since $\mathbb{T}^n$ is commutative, the brackets $[X_i,X_j]$ vanish for all $i$ and $j$. Hence for any $\lambda _{12}$ in $\mathbb{R}$, the twist $t_{12}:=\lambda _{12}X_1\wedge X_2$ deforms $\omega _{\mathrm{FS}}$ to a symplectic structure $\omega _{\mathrm{FS}}^{t_{12}}$ induced by a Poisson structure $\pi _{\mathrm{FS}}^{t_{12}}:=\pi _{\mathrm{FS}}-\left(t_{12}\right)_{\mathbb{C}\mathrm{P}^n}$ by Theorem \ref{main}. On the other hand it follows $\pi _{\mathbb{T}^n}^t=t^L-t^R=0$ for any twist $t$ by the commutativity of $\mathbb{T}^n$. Therefore, after deformation, the multiplicative Poisson structure $0$ on $\mathbb{T}^n$ is invariant and this action is a symplectic-Hamiltonian action with the same moment map $\mu $. Therefore, by Theorem \ref{main} again, the twist $t_{13}:=\lambda _{13}X_1\wedge X_3$ deforms $\omega _{\mathrm{FS}}^{t_{12}}$ to $(\omega _{\mathrm{FS}}^{t_{12}})^{t_{13}}=\omega _{\mathrm{FS}}^{t_{12}+t_{13}}$ induced by $(\pi _{\mathrm{FS}}^{t_{12}})^{t_{13}}
%=\pi _{\mathrm{FS}}^{t_{12}}-\left(t_{13}\right)_{\mathbb{C}\mathrm{P}^n}=\pi _{\mathrm{FS}}-\left(t_{12}+t_{13}\right)_{\mathbb{C}\mathrm{P}^n}
=\pi _{\mathrm{FS}}^{t_{12}+t_{13}}$. Then we see that the trivial Poisson structure on $\mathbb{T}^n$ is invariant and that the action is a symplectic-Hamiltonian action with $\mu $. By repeating this operation, it follows that we can deform $\omega _{\mathrm{FS}}$ to $\omega _{\mathrm{FS}}^t$ for any twist $t=\sum _{i<j}^{}\lambda _{ij}X_i\wedge X_j$ and that the action is a symplectic-Hamiltonian action with $\mu $. On $U_1$, since we obtain
\begin{eqnarray*}
(X_i\wedge X_j)_{\mathbb{C}\mathrm{P}^n}=y_iy_j\frac{\partial}{\partial x_i}\wedge \frac{\partial}{\partial x_j}-y_ix_j\frac{\partial}{\partial x_i}\wedge \frac{\partial}{\partial y_j}-x_iy_j\frac{\partial}{\partial y_i}\wedge \frac{\partial}{\partial x_j}+x_ix_j\frac{\partial}{\partial y_i}\wedge \frac{\partial}{\partial y_j}\\
(1\leq i<j\leq n),
\end{eqnarray*}
where $x_i:=\mathrm{Re}\frac{z_{i+1}}{z_1}$ and $y_i:=\mathrm{Im}\frac{z_{i+1}}{z_1}$, it follows that 
\begin{eqnarray*}
\omega_{\mathrm{FS}}^t&=&\omega_{\mathrm{FS}}+\sum_{i<j}^{}\frac{\lambda _{ij}}{\left\{\sum_{k}^{}(x_k^2+y_k^2)+1\right\}^3}(x_ix_jdx_i\wedge dx_j\\
              &&\qquad \qquad \qquad \qquad \qquad \quad +x_iy_jdx_i\wedge dy_j
+y_ix_jdy_i\wedge dx_j+y_iy_jdy_i\wedge dy_j).
\end{eqnarray*}
\end{example}

The last example is the complex Grassmannian $\mathrm{Gr}(r;\mathbb{C}^n):=\mathrm{SU}(n)/(\mathrm{S}(\mathrm{U}(r)\times \mathrm{U}(n-r)))$ with the Kirillov-Kostant form $\omega _{\mathrm{KK}}$. %Here $\omega _{\mathrm{KK}}$ is a symplectic form naturally defined on a coadjoint orbit of $\mathrm{SU}(n)$ and  $\mathrm{Gr}(r;\mathbb{C}^n)$ with the natural action is identified with one of coadjoint orbits of $\mathrm{SU}(n)$ with the coadjoint action. 
With respect to $\omega _{\mathrm{KK}}$, the natural $\mathrm{SU}(n)$-action is symplectic-Hamiltonian (\cite{Si}).

Then we consider the following r-matrix of $\mathfrak{su}(n)$:
$$t=\frac{1}{4n}\sum _{1\leq i<j\leq n}^{}X_{ij}\wedge Y_{ij},$$
where the r-matrix $t$ is the canonical one defined on any compact semi-simple Lie algebra over $\mathbb{R}$ (for example, see \cite{D}).
This is an r-matrix such that $[t,t]\neq 0$. We show that it satisfies $[t,t]_M=0$, where $M:=\mathrm{Gr}(r;\mathbb{C}^n)$. Since $t$ is an r-matrix, the element $[t,t]$ is $\mathrm{ad}$-invariant by the definition. Therefore $[t,t]$ is $\mathrm{Ad}$-invariant because $\mathrm{SU}(n)$ is connected. By the definition of the $\mathrm{SU}(n)$-action on $\mathrm{Gr}(r;\mathbb{C}^n)$, it follows that
\begin{eqnarray*}
[t,t]_M=p_*[t,t]^R,
\end{eqnarray*}
where $p:\mathrm{SU}(n)\rightarrow \mathrm{Gr}(r;\mathbb{C}^n)=\mathrm{SU}(n)/(\mathrm{S}(\mathrm{U}(r)\times \mathrm{U}(n-r)))$ is the natural projection. Since any point $m$ in $\mathrm{Gr}(r;\mathbb{C}^n)$ is represented by $gH$, where $g$ in $\mathrm{SU}(n)$ and $H:=\mathrm{S}(\mathrm{U}(r)\times \mathrm{U}(n-r))$, we compute
\begin{eqnarray*}
[t,t]_{M,m}=p_*[t,t]_g^R=p_*R_{g*}[t,t].
\end{eqnarray*}
Because of the $\mathrm{Ad}$-invariance of $[t,t]$, we obtain
\begin{eqnarray*}
p_*R_{g*}[t,t]=p_*L_{g*}L_{g^{-1}*}R_{g*}[t,t]=p_*L_{g*}\mathrm{Ad}_{g^{-1}}[t,t]=p_*L_{g*}[t,t].
\end{eqnarray*}
Let $\mathfrak{h}$ be the Lie algebra of $H$. For any $X$ in $\mathfrak{h}$ and $g$ in $\mathrm{SU}(n)$, we compute
\begin{eqnarray*}
p_*L_{g*}X=p_*L_{g*}\left.\frac{d}{ds}\exp sX\right|_{s=0}
          =\left.\frac{d}{ds}(g\exp sX)H\right|_{s=0}
          =\left.\frac{d}{ds}gH\right|_{s=0}
          =0,
\end{eqnarray*}
where we have used that $\exp sX$ is in $H$ in the third equality. Therefore it holds that $[t,t]_M=0$ if each term of $[t,t]$ includes elements in $\mathfrak{h}$ as follows. We notice that
$$\mathfrak{h}=\mathrm{span}_\mathbb{R}\{X_{ij},Y_{ij},Z_k|1\leq i<j\leq r \ \mathrm{or}\ r+1\leq i<j\leq n,\ \mathrm{and}\ k=1,\dots ,n-1\}.$$
If $X_{ij},Y_{ij}\in \mathfrak{h}$, then
$$\left[\ \cdot \ ,X_{ij}\wedge Y_{ij}\right]=\left[\ \cdot \ ,X_{ij}\right]\wedge Y_{ij}-X_{ij}\wedge \left[\ \cdot \ ,Y_{ij}\right].$$
So these terms include an element in $\mathfrak{h}$. Hence we investigate terms of the form
\begin{eqnarray*}
\left[X_{ij}\wedge Y_{ij},X_{kl}\wedge Y_{kl}\right]=&-\left[X_{ij},X_{kl}\right]\wedge Y_{ij}\wedge Y_{kl}-X_{ij}\wedge \left[Y_{ij},X_{kl}\right]\wedge Y_{kl}\\
                  &\quad -Y_{ij}\wedge \left[X_{ij},Y_{kl}\right]\wedge X_{kl}-X_{ij}\wedge X_{kl}\wedge \left[Y_{ij},Y_{kl}\right],
\end{eqnarray*}
where $X_{ij}, Y_{ij},X_{kl}$ and $Y_{kl}$ are not in $\mathfrak{h}$. In the case of $i=k$ and $j=l$, we get
\begin{eqnarray*}
\left[X_{ij},X_{ij}\right]=\left[Y_{ij},Y_{kl}\right]=0,\\
\left[X_{ij},Y_{ij}\right]=2(Z_i-Z_j)\in \mathfrak{h},
\end{eqnarray*}
where $Z_n=0$. In the case of  $i=k$ and $j<l$ (resp. $l<j$), since it follows that $r+1\leq j,\ l\leq n$, we obtain
\begin{eqnarray*}
\left[X_{ij},X_{kl}\right]=\left[Y_{ij},Y_{kl}\right]=-X_{jl}(\mathrm{resp.}\ X_{lj})\in \mathfrak{h},\\
\left[Y_{ij},X_{kl}\right]=\left[Y_{kl},X_{ij}\right]=-Y_{jl}(\mathrm{resp.}\ Y_{lj})\in \mathfrak{h}.
\end{eqnarray*}
We can also show the case of $i<k$ (resp. $k<i$) and $j=l$ in the similar way.
%%%%%%%%%%%%%%%%%%%%%%%%%%%%%%%%%%%%%%%%%%%%%%%%%%%%%%%%%%%%%%%%%%%%%%
\iffalse
we compute
\begin{eqnarray*}
\left[X_{ij},X_{kl}\right]=\left[Y_{ij},Y_{kl}\right]=-X_{\alpha _{ik}}(\mathrm{resp.}\ X_{\alpha _{ki}}),\\
\left[Y_{ij},X_{kl}\right]=\left[Y_{kl},X_{ij}\right]=Y_{\alpha _{ik}}(\mathrm{resp.}-Y_{\alpha _{ki}}).
\end{eqnarray*}
These are in $\mathfrak{h}$ since it follows that $1\leq i,\ k\leq r$. 
\fi
%%%%%%%%%%%%%%%%%%%%%%%%%%%%%%%%%%%%%%%%%%%%%%%%%%%%%%%%%%%%%%%%%%%%%%%%
Therefore all terms of $[t,t]$ include elements in $\mathfrak{h}$, so that $[t,t]_M=0$. Since $\mathrm{Gr}(r;\mathbb{C}^n)$ is compact, for sufficiently small $|\lambda |$, a $2$-vector field $\pi _{\mathrm{KK}}^{\lambda t}$ is nondegenerate by Theorem \ref{abstract main}, where $\pi _{\mathrm{KK}}$ is the Poisson structure induced by $\omega _{\mathrm{KK}}$. Example \ref{cp1} is the special case of this example. From the above discussion, we obtain the following.

\begin{theorem}
Let $t$ be the above r-matrix of $\mathfrak{su}(n)$. Then there exists sufficiently small number $\lambda $ such that the Kirillov-Kostant form $\omega _{\mathrm{KK}}$ on $\mathrm{Gr}(r;\mathbb{C}^n)$ can be deformed by a twist $\lambda t$ in the sense of Section \ref{Main Result}.
\end{theorem}

\section*{acknowledgments}
I would like to express my deepest gratitude to Yuji Hirota for leading me into the study of quasi-Poisson theory and my supervisor Yasushi Homma for his helpful advice.
%\end{acknowledgments}

%\section*{References}


\begin{thebibliography}{99}
\bibitem{AK}
{A. Alekseev and Y. Kosmann-Schwarzbach}.
Manin pairs and moment maps. 
\emph{J. Diff. Geom. }
\textbf{56} (2000) 133--165.

\bibitem{BC}
{H. Bursztyn and M. Crainic}. 
Dirac geometry, quasi-Poisson actions and $D/G$-valued moment maps. 
\emph{J. Diff. Geom. }
\textbf{82} 3 (2009) 501--566.

\bibitem{D}
{V. G. Drinfel'd}. 
Hamiltonian structures on Lie groups, Lie bialgebras and the geometric meaning of the classical Yang-Baxter equations. 
\emph{Soviet Math. Dokl. }
\textbf{27} (1983) 68--71.

\bibitem{L}
{J. -H. Lu}. 
Momentum mappings and reduction of Poisson actions, in Symplectic Geometry, Groupoids and Integrable Systems, P. Dazord and A. Weinstein, eds.  
\emph{Springer} (1991) 209--226.

\bibitem{L2}
{J. -H. Lu}. 
Multiplicative and Affine Poisson structures on Lie groups.  
Berkeley Thesis (1991).

\bibitem{LW}
{J. -H. Lu and A. Weinstein}. 
Poisson-Lie groups, dressing transformations and Bruhat decompositions. 
\emph{J. Diff. Geom. }
\textbf{31} (1990) 501--526.

\bibitem{MR}
{J. E. Marsden and T. S. Ratiu}. 
\emph{Introduction to Mechanics and Symmetry Second Edition.}
(Springer, 2003).

\bibitem{Sa}
{D. Salamon}. 
Uniqueness of Symplectic Structures. 
arXiv:1211.2940v5

\bibitem{Si}
{A. C. da Silva}. 
\emph{Lectures on Symplectic Geometry.}  
(Springer-Verlag, 2006).

\bibitem{V}
{I. Vaisman}. 
Lectures on the Geometry of Poisson Manifold.  
\emph{Birkhaeuser} (1994).

\end{thebibliography}
\end{document}